\documentclass[12pt,leqno]{amsart} 
\newtheorem{theo}{Theorem}[section]
\newtheorem{defi}{Definition}[section]
\newtheorem{lem}{Lemma}[section]
\newtheorem{cor}{Corollary}[section]
\newtheorem{rem}{Remark}[section]
\newcommand{\be}{\begin{equation}}
\newcommand{\ee}{\end{equation}}
\newcommand{\bea}{\begin{eqnarray}}
\newcommand{\eea}{\end{eqnarray}}
\usepackage{amssymb,amsfonts,amsthm,indentfirst}
\numberwithin{equation}{section}
\begin{document}
\title { ON LOCALLY $\phi$-SEMISYMMETRIC SASAKIAN MANIFOLDS} 
\author[A. A. Shaikh, C. K. Mondal and H. Ahmad]{Absos Ali Shaikh, Chandan Kumar Mondal and Helaluddin Ahmad }
\date{}
\address{\noindent\newline Department of Mathematics,\newline University of 
Burdwan, Golapbag,\newline Burdwan-713104,\newline West Bengal, India}
\email{\small{aask2003@yahoo.co.in, aashaikh@math.buruniv.ac.in}}

\noindent\footnotetext {$\mathbf{2000}$\hspace{5pt}Mathematics\; Subject\; Classification,\; 53C15, 53C25. }
\begin{abstract} 
Generalizing the notion of local $\phi$-symmetry of Takahashi \cite{TAKA}, in the present paper, we introduce the notion of  \textit{local $\phi$-semisymmetry} of a Sasakian manifold along with its proper existence and characterization. We also study the notion of  local Ricci (resp., projective, conformal) $\phi$-semisymmetry of a  Sasakian manifold and obtain its characterization. It is shown that the local $\phi$-semisymmetry, local projective $\phi$-semisymmetry and local concircular $\phi$-semisymmetry are equivalent. It is also shown that local conformal $\phi$-semisymmetry and local conharmonical $\phi$-semisymmetry are equivalent.
\end{abstract} 
\thanks{Typeset by \AmS-\LaTeX}
\maketitle
\flushbottom

 \section {\bf Introduction} 
\indent Let $M$ be an $n$-dimensional, $n \geq 3$, connected smooth Riemannian manifold endowed with the Riemannian metric $g$. Let $\nabla$, $R$, $S$ and $r$ be the  Levi-Civita connection, curvature tensor, Ricci tensor and the scalar curvature of $M$ respectively. The manifold $M$ is called locally symmetric due to Cartan (\cite{cartan}, \cite{cartan1}) if the local geodesic symmetry at $p\in M$ is an isometry, which is equivalent to the fact that $\nabla R=0$. Generalizing the concept of local symmetry, the notion of semisymmetry was introduced by Cartan \cite{cartan2} and fully classified by Szab\'o (\cite{szabo1}, \cite{szabo2}, \cite{szabo3}). The manifold $M$ is said to be semisymmetric if
\be\nonumber
(R(U,V).R)(X,Y)Z = 0
\ee
for all vector fields $X$, $Y$, $Z$, $U$,  $V$ on $M$, where $R(U,V)$ is considered as the derivation of the tensor algebra at each point of $M$. Every locally symmetric manifold is semisymmetric but the converse is not true, in general. However, the converse is true only for $n=3$.
  As a weaker version of local symmetry, in 1977 Takahashi \cite{TAKA} introduced the notion of local $\phi$-symmetry on a Sasakian manifold. A Sasakian manifold  is said to be locally $\phi$-symmetric  if
\be\nonumber
\phi^2((\nabla_WR)(X,Y)Z) = 0
\ee
for all horizontal vector fields $X$, $Y$, $Z$, $W$ on $M$, where $\phi$ is the structure tensor of the manifold $M$. The concept of local $\phi$-symmetry on various structures and their generalizations or extensions are studied in \cite{ds}, \cite{sb}, \cite{sbe1}, \cite{sbe2}, \cite{ste1}, \cite{ste2},  \cite{sd}, \cite{sh1}, \cite{sh2}. By extending the notion of semisymmetry and generalizing the concept of local $\phi$-symmetry of Takahashi \cite{TAKA}, in the present paper, we introduce the notion of \textit{local $\phi$-semisymmetry} on a Sasakian manifold. A Sasakian manifold $M$, $n\geq 3$, is said to be \textit{locally $\phi$-semisymmetric} if
\be\nonumber
\phi^2((R(U,V).R)(X,Y)Z) = 0
\ee
for all horizontal vector fields $X$, $Y$, $Z$, $U$,  $V$ on $M$. We note that every locally $\phi$-symmetric as well as semisymmetric Sasakian manifold is locally $\phi$-semisymmetric but not conversely. The object of the present paper is to study the geometric properties of a locally $\phi$-semisymmetric Sasakian manifold  along with its proper existence and characterization. The paper is organized as follows. Section~2 deals with the rudiments of Sasakian manifolds. By extending the definition of local $\phi$-symmetry, in Section~3, we derive the defining condition of a locally $\phi$-semisymmetric Sasakian manifold and proved that a Sasakian manifold is locally $\phi$-semisymmetric if and only if each K\"ahlerian manifold, which is a base space of a local fibering, is Hermitian locally semisymmetric. We cite an example of a locally $\phi$-semisymmetric Sasakian manifold which is not locally $\phi$-symmetric. We also obtain a characterization of locally $\phi$-semisymmetric Sasakian manifold by considering the horizontal vector fields. Section~4 is devoted to the characterization of locally $\phi$-semisymmetric Sasakian manifold for arbitrary vector fields. As the generalization of Ricci (resp., projectively, conformally)  semisymmetric Sasakian manifold, in the last section, we introduce the notion of  \textit{locally Ricci  (resp., projectively, conformally) $\phi$-semisymmetric Sasakian manifold} and obtain the characterization of such notions. Recently Shaikh and Kundu \cite{kundu} defined a generalized curvature tensor, called $B$-tensor, by the linear combination of $R$, $S$ and $g$ which includes various curvature tensors as particular cases. We study the characterization of locally $B$-$\phi$-semisymmetric Sasakian manifolds. It is shown that local $\phi$-semisymmetry, local projective $\phi$-semisymmetry and local concircular $\phi$-semisymmetry are equivalent and hence they are of the same characterization. Also it is  proved that local conformal $\phi$-semisymmetry and local conharmonical $\phi$-semisymmetry are equivalent. Finally, we conclude that the study of local $\phi$-semisymmetry and local conformal $\phi$-semisymmetry are meaningful as they are not equivalent. However, the study of local $\phi$-semisymmetry with any other  generalized curvature tensor of type (1,3) (which are the linear combination of $R$, $S$ and $g$) is either meaningless or redundant due to their equivalency.  
\section{Sasakian manifolds}
\indent An $n (= 2m+1, m\geq 1)$-dimensional $C^{\infty}$ manifold $M$ is said to be a contact manifold if it carries a global 1-form $\eta$ such that $\eta \wedge (d\eta)^{m} \neq 0$ everywhere on $M$. Given a contact form $\eta$, it is well-known that there exists a unique vector field $\xi$, called the characterstic vector field of $\eta$, satisfying $\eta(\xi)=1$ and $d\eta(X,\xi)=0$ for any vector field $X$ on $M$. A Riemannian metric $g$ is said to be an associated metric if there exists a tensor field $\phi$ of type (1,1) such that   
\be
\label{eq:2.1}
\phi ^{2} = -I+\eta \otimes\xi, \ \ \eta(\cdot) = g(\cdot, \xi), \ \ \ d\eta(\cdot,\cdot)=g(\cdot,\phi \cdot).
\ee
Then the structure $(\phi, \xi, \eta, g)$ on $M$ is called a contact metric stucture and the manifold $M$  equipped with such a stucture is called a contact metric manifold~\cite{BLAIR}.\\
From (\ref{eq:2.1}) it is easy to check that the following holds:
\begin{eqnarray}
\label{eq:2.2} 
&&\phi\xi = 0, \ \ \ \eta\circ\phi= 0, \ \ \  g(\phi \cdot,\cdot) = -g(\cdot,\phi \cdot),\\ 
\label{eq:2.3}
&&g(\phi \cdot, \phi \cdot) = g(\cdot, \cdot) - \eta \otimes\eta .
\end{eqnarray}
\indent Given a contact metric manifold $M$ there is an (1,1) tensor field $h$ given by $h=\frac{1}{2}\pounds_\xi\phi$, where $\pounds$ denotes the operator of Lie differentiation. Then $h$ is symmetric. The vector field $\xi$ is a Killing vector field with respect to $g$ if and only if $h=0$. A contact metric manifold $M$  for which $\xi$ is a Killing vector is said to be a $K$-contact manifold. A contact structure on $M$ gives rise to an almost complex structure $J$ on the product $M\times \mathbb R$ defined by 
\begin{equation*}
J\Big{(}X,f\frac{d}{dt}\Big{)}=\Big{(}\phi X-f\xi,\eta(X)\frac{d}{dt}\Big{)}, 
\end{equation*}
where $f$ is a real valued function, is integrable, then the structure is said to be normal and the manifold $M$ is a Sasakian manifold. Equivalently, a contact metric manifold is Sasakian if and only if
\begin{equation}
\label{eq:2.4}
R(X,Y)\xi = \eta(Y)X-\eta(X)Y
\end{equation}
holds for all $X$, $Y$ on $M$.\\
In an $n$-dimensional Sasakian manifold $M$ the following relations hold (\cite{BLAIR}, \cite{KON}):
\bea
\label{eq:2.5}
&&R(\xi,X)Y = (\nabla _{X}\phi)(Y) = g(X,Y)\xi - \eta(Y)X = - R(X,\xi)Y,\\
\label{eq:2.6}
&&\nabla_{X}\xi = -\phi X, \ \ \ (\nabla _{X}\eta)(Y) = g(X,\phi Y),\\
\label{eq:2.7}
&&\eta(R(X,Y)Z) = g(Y,Z)\eta(X) - g(X,Z)\eta(Y),\\
\label{eq:2.8}
&&(\nabla_WR)(X,Y)\xi = g(W,\phi Y)X - g(W,\phi X)Y + R(X,Y)\phi W,
\eea
\bea
\label{eq:2.9}
&&(\nabla_WR)(X,\xi)Z = g(X,Z)\phi W - g(Z,\phi W)X + R(X,\phi W)Z,\\
\label{eq:2.10}
&&S(X ,\xi) = (n-1) \eta(X), \ \ S(\xi, \xi) = (n-1)
\eea
for all vector fields $X$, $Y$, $Z$ and $W$ on $M$. In a Sasakian manifold, for any $X$, $Y$, $Z$ on $M$, we also have \cite{tanno}
\be\label{eq:2.11}
R(X,Y)\phi W= g(W,\phi X)Y - g(W,Y)\phi X - g(W,\phi Y)X + g(W,X)\phi Y +\phi R(X,Y)W.
\ee
From (\ref{eq:2.8}) and (\ref{eq:2.11}), it follows that
\be\label{eq:2.12}
(\nabla_WR)(X,Y)\xi =g(W,X)\phi Y  - g(W,Y)\phi X + \phi R(X,Y)W.             
\ee
\section{\bf Locally $\phi$-Semisymmetric Sasakian Manifolds}
Let $M$ be an $n (= 2m+1, m\geq 1)$-dimensional Sasakian manifold endowed with the structure $(\phi, \xi, \eta, g)$. Let $\tilde U$ be an  open neighbourhood of $x \in M$ such that the induced Sasakian structure on $\tilde U$, denoted by the same letters, is regular. Let $\pi : \tilde U\rightarrow N=\tilde U/\xi$ be a (local) fibering and let $(J, \bar{g})$ be the induced K\"ahlerian structure on $N$ \cite{og}. Let $R$ and $\bar{R}$ be the curvature tensors constructed by $g$ and $\bar{g}$ respectively. For a vector field $\bar{X}$ on $N$, we denote its horizontal lift (with respect to the connection form $\eta$) by $\bar{X}^*$. Then we have, for any vector fields $\bar{X}$, $\bar{Y}$ and $\bar{Z}$ on $N$,
\be\label{eq:3.1}
(\bar{\nabla}_{\bar{X}}{\bar{Y}})^{*} = \nabla_{{\bar{X}}^{*}}{\bar{Y}}^{*} - \eta(\nabla_{{\bar{X}}^{*}}{\bar{Y}}^{*})\xi,
\ee   
\be\label{eq:3.2}
(\bar{R}(\bar{X},\bar{Y})\bar{Z})^{*}=R(\bar{X}^{*},\bar{Y}^{*})\bar{Z}^{*}+g(\phi {\bar{Y}^{*}},\bar{Z}^{*})\phi {\bar{X}^{*}} - g(\phi {\bar{X}^{*}},\bar{Z}^{*})\phi {\bar{Y}^{*}}-2g(\phi {\bar{X}^{*}},\bar{Y}^{*})\phi {\bar{Z}^{*}},
\ee
\be\label{eq:3.3}
(({\bar{\nabla}}_{\bar{V}} \bar{R})(\bar{X}, \bar{Y})\bar{Z})^{*}=-\phi^{2} [(\nabla_{\bar{V}^{*}} R)(\bar{X}^{*}, \bar{Y}^{*})\bar{Z}^{*}]
\ee 
where $\bar{\nabla}$ is the Levi-Civita connection for $\bar{g}$. The relations (\ref{eq:3.1}) and (\ref{eq:3.2}) are due to Ogiue \cite{og} and the relation (\ref{eq:3.3}) is due to Takahashi \cite{TAKA}.\\
\indent Making use of (\ref{eq:2.1}), (\ref{eq:2.4})-(\ref{eq:2.11}) and (\ref{eq:3.1})-(\ref{eq:3.3}), we get by straightforward calculation
\be\label{eq:3.4}
((\bar{R}(\bar{U},\bar{V})\cdot \bar{R})(\bar{X},\bar{Y})\bar{Z})^{*}=-\phi^{2} [(R(\bar{U}^{*},\bar{V}^{*})\cdot R)(\bar{X}^{*},\bar{Y}^{*})\bar{Z}^{*}]
\ee
for any vector fields $\bar{X},\bar{Y}, \bar{Z}, \bar{U}$ and $\bar{V}$ on $N$, where $R(U,V)$ is considered as the derivation of the tensor algebra at each point of $N$. Hence from (\ref{eq:3.4}) it is natural to define the following:
\begin{defi}
A Sasakian manifold is said to be a locally $\phi$-semisymmetric if
\be\label{eq:3.5}
\phi^{2}[(R(U,V)\cdot R)(X,Y)Z]=0
\ee
for any horizontal vector fields $X, Y, Z, U$ and $V$ on $M$, where a horizontal vector is a vector which is horizontal with respect to the connection form $\eta$ of the local fibering, that is, orthogonal to $\xi$.
\end{defi}
Thus from (\ref{eq:3.4}) and (\ref{eq:3.5}), we can state the following:
\begin{theo}
A Sasakian manifold is locally $\phi$-semisymmetric if and only if each K\"ahlerian manifold, which is a base space of a local fibering, is a Hermitian locally semisymmetric space.
\end{theo}
\noindent{\bf Example 3.1.} Suppose a Sasakian manifold is not of constant $\phi$-sectional curvature. Then the K\"ahlerian base manifold is not of constant sectional curvature.\\
Now let $R(X,\phi X,Y,\phi Y)=f\in C^{\infty}(M)$. Then $(\nabla_V R)(X,\phi X,Y,\phi Y) = (Vf) \ne 0$, i.e. the K\"ahlerian manifold is not Hermitian locally symmetric and therefore the Sasakian manifold is not locally $\phi$-symmetric.\\
Now $(\nabla_U \nabla_V R)(X,\phi X,Y,\phi Y) = U(Vf) \ne 0$, which implies that $(R(U,V)\cdot R)(X,\phi X,Y,\phi Y)=0$,
i.e. the K\"ahlerian manifold is Hermitian locally semisymmetric.\\
Hence the Sasakian manifold is locally $\phi$-semisymmetric but not locally $\phi$-symmetric.\\
\indent First we suppose that $M$ is a Sasakian manifold such that
\be\label{eq:3.6}
\phi^{2}[(R(U,V)\cdot R)(X,Y)\xi]=0
\ee
for any horizontal vector fields $X, Y, U$ and $V$ on $M$.\\
Differentiating (\ref{eq:2.12}) covariantly with respect to a horizontal vector field $U$ we get
\bea\label{eq:3.7}
(\nabla_U\nabla_VR)(X,Y)\xi &=& \{g(Y,U)g(X,V) - g(X,U)g(Y,V) - R(X,Y,U,V)\}\xi\\
\nonumber &+&\phi ((\nabla_UR)(X,Y)V).             
\eea
Alternating $U$ and $V$ on (\ref{eq:3.7}) we get
\bea\label{eq:3.8}
(\nabla_V\nabla_UR)(X,Y)\xi &=& \{g(Y,V)g(X,U) - g(X,V)g(Y,U) - R(X,Y,V,U)\}\xi\\
\nonumber &+&\phi ((\nabla_VR)(X,Y)U).             
\eea
From  (\ref{eq:3.7}) and (\ref{eq:3.8}) it follows that  
\bea\label{eq:3.9}
(R(U,V)\cdot R)(X,Y)\xi&=&2\{g(Y,U)g(X,V) - g(X,U)g(Y,V) - R(X,Y,U,V)\}\xi\\
\nonumber &+& \phi \{(\nabla_UR)(X,Y)V-(\nabla_VR)(X,Y)U\}.
\eea
Again from (\ref{eq:3.6}) we have
\be\label{eq:3.10}
(R(U,V)\cdot R)(X,Y)\xi=0.
\ee
From (\ref{eq:3.9}) and (\ref{eq:3.10}) we have
\bea\label{eq:3.11}
2\{g(Y,U)g(X,V) &-& g(X,U)g(Y,V) - R(X,Y,U,V)\}\xi\\
\nonumber&+&\phi \{(\nabla_UR)(X,Y)V-(\nabla_VR)(X,Y)U\}=0.
\eea
Applying $\phi$ on (\ref{eq:3.11}) and using (\ref{eq:2.11}), (\ref{eq:2.12}) and (\ref{eq:2.2}) we get
\be\label{eq:3.12}
(\nabla_UR)(X,Y)V-(\nabla_VR)(X,Y)U=0.
\ee
In view of (\ref{eq:3.12}), (\ref{eq:3.11}) yields
\be\label{eq:3.13}
 R(X,Y,U,V)=g(Y,U)g(X,V) - g(X,U)g(Y,V)
\ee
for any horizontal vector fields $X, Y, U$ and $V$ on $M$. Hence $M$ is of constant $\phi$-holomorphic sectional curvature 1 and hence of constant curvature 1. This leads to the following:
\begin{theo}
If a Sasakian manifold $M$ satisfies the condition $\phi^{2}[(R(U,V)\cdot R)(X,Y)\xi]=0$ for all horizontal vector fields $X, Y, Z, U$ and $V$  on $M$, then it is a  manifold of constant curvature 1.
\end{theo}
Now we consider a locally $\phi$-semisymmetric Sasakian manifold. Then from (\ref{eq:3.5}) we have
\be\nonumber
(R(U,V)\cdot R)(X,Y)Z=g((R(U,V)\cdot R)(X,Y)Z,\xi)\xi,
\ee
from which we get
\be\label{eq:3.14}
(R(U,V)\cdot R)(X,Y)Z=-g((R(U,V)\cdot R)(X,Y)\xi,Z)\xi
\ee
for all horizontal vector fields $X, Y, Z, U$ and $V$ on $M$.\\
In view of (\ref{eq:3.9}), it follows from (\ref{eq:3.14}) that
\be\label{eq:3.15}
(R(U,V)\cdot R)(X,Y)Z=[(\nabla_UR)(X,Y,V,\phi Z)-(\nabla_VR)(X,Y,U,\phi Z)]\xi.
\ee
Now differentiating (\ref{eq:2.11}) covariantly with respect to a horizontal vector field $V$,  we obtain
\bea\label{eq:3.16}
(\nabla_VR)(X,Y)\phi Z&=&[R(X,Y,Z,V)-\{g(Y,Z)g(X,V) + g(X,Z)g(Y,V)]\xi \\
\nonumber &+&\phi ((\nabla_VR)(X,Y)Z). 
\eea
Taking inner product of (\ref{eq:3.16}) with a horizontal vector field $U$, we obtain
\be\label{eq:3.17}
g((\nabla_VR)(X,Y)\phi Z,U)=-g((\nabla_VR)(X,Y)Z,\phi U).
\ee
Using (\ref{eq:3.17}) in (\ref{eq:3.15}) we get
\be\label{eq:3.18}
(R(U,V)\cdot R)(X,Y)Z=[(\nabla_UR)(X,Y,Z,\phi V )-(\nabla_VR)(X,Y,Z,\phi U)]\xi
\ee
for any horizontal vector fields $X, Y, Z, U$ and $V$ on $M$. Hence we can state the following:
\begin{theo}
A necessary and sufficient condition for a Sasakian manifold $M$ to be a locally $\phi$-semisymmetric is that it satisfies the relation \textnormal{(\ref{eq:3.18})} for all horizontal vector fields on $M$.
\end{theo}
\section{\bf Characterization of  a Locally $\phi$-Semisymmetric Sasakian Manifold}
In this section we investigate the condition of local $\phi$-semisymmetry of a Sasakian manifold for arbitrary vector fields on $M$. To find the condition we need the following lemmas.
\begin{lem}\cite{TAKA}
For any horizontal vector fields $X, Y$ and $Z$ on $M$, we get 
\be\label{eq:4.1}
(\nabla_\xi R)(X,Y)Z=0.
\ee
\end{lem}
Now  Lemma 4.1, (\ref{eq:2.9}) and (\ref{eq:2.12}) together imply the following:
\begin{lem}\cite{TAKA}
For any  vector fields $X, Y, Z, V$ on $M$, we get
\bea\label{eq:4.2}
(\nabla_{\phi^2 V}R)(\phi^2 X,\phi^2 Y)\phi^2 Z&=&(\nabla_VR)(X,Y)Z \\ 
\nonumber &+& \eta(X)\{g(Y,Z)\phi V-g(\phi V,Z)Y+R(Y,\phi V)Z\}\\
\nonumber &-& \eta(Y)\{g(X,Z)\phi V-g(\phi V,Z)X + R(X,\phi V)Z\}\\
\nonumber &-& \eta(Z)\{g(X,V)\phi Y-g(Y,V)\phi X + \phi R(X,Y)V\}.
\eea
\end{lem}
\indent Now let $X, Y, Z, U, V$ be arbitrary vector fields on $M$. We now compute  $(R(\phi^2 U,\phi^2 V)\cdot R)(\phi^2 X,\phi^2 Y)\phi^2 Z$ in two different ways. Firstly, from (\ref{eq:3.18}), (\ref{eq:2.1}) and (\ref{eq:4.2}) we get
\bea\label{eq:4.3} 
&&(R(\phi^2 U,\phi^2 V)\cdot R)(\phi^2 X,\phi^2 Y)\phi^2 Z=[\{(\nabla_VR)(X,Y,Z,\phi U) - (\nabla_UR)(X,Y,Z,\phi V)\}\\
\nonumber &&+\eta(X)\{ g(Y,\phi V)g(\phi U,Z)- g(Y,\phi U)g(\phi V,Z) - R(Y,Z,\phi U,\phi V)\}\\
\nonumber &&-\eta(Y)\{ g(X,\phi V)g(\phi U,Z)- g(X,\phi U)g(\phi V,Z) - R(X,Z,\phi U,\phi V)\}\\
\nonumber &&-2\eta(Z)\{ g(X,V)g(U,Y)- g(X,U)g(V,Y) - R(X,Y,U,V)\}]\xi.
\eea
Again using  (\ref{eq:2.11}) in (\ref{eq:4.3}), we obtain  
\bea\label{eq:4.4} 
&&(R(\phi^2 U,\phi^2 V)\cdot R)(\phi^2 X,\phi^2 Y)\phi^2 Z=[\{(\nabla_VR)(X,Y,Z,\phi U) - (\nabla_UR)(X,Y,Z,\phi V)\}\\
\nonumber &&-  \eta(X)H(Y,Z,U,V)+ \eta(Y)H(X,Z,U,V)+ 2\eta(Z)H(X,Y,U,V)]\xi
\eea
where $H(X,Y,Z,U)=g(\mathcal{H}(X,Y)Z,U)$ and the tensor field $\mathcal{H}$ of type (1,3) is given by
\be\label{eq:4.5}
\mathcal{H}(X,Y)Z=R(X,Y)Z-g(Y,Z)X+g(X,Z)Y
\ee
for all vector fields $X, Y, Z$ on $M$.
Secondly, we have 
\bea\label{eq:4.6} 
&&(R(\phi^2 U,\phi^2 V)\cdot R)(\phi^2 X,\phi^2 Y)\phi^2 Z=R(\phi^2 U,\phi^2 V)R(\phi^2 X,\phi^2 Y)\phi^2 Z\\
\nonumber&&-R(R(\phi^2 U,\phi^2 V)\phi^2 X,\phi^2 Y)\phi^2 Z- R(\phi^2 X,R(\phi^2 U,\phi^2 V)\phi^2 Y)\phi^2 Z\\
\nonumber&&-R(\phi^2 X,\phi^2 Y)R(\phi^2 U,\phi^2 V)\phi^2 Z.
\eea
By straightforward calculation, from (\ref{eq:4.6}) we get 
\bea\label{eq:4.7} 
(R(\phi^2 U,\phi^2 V)\cdot R)(\phi^2 X,\phi^2 Y)\phi^2 Z &=& - (R(U,V)\cdot R)(X,Y)Z\\
\nonumber&+& \eta(U)\big{[}H(X,Y,Z,V)\xi+\eta(X)\mathcal{H}(V,Y)Z\\
\nonumber&+& \eta(Y)\mathcal{H}(X,V)Z+\eta(Z)\mathcal{H}(X,Y)V\big{]}\\
\nonumber&-& \eta(V)\big{[}H(X,Y,Z,U)\xi+\eta(X)\mathcal{H}(U,Y)Z\\
\nonumber&+& \eta(Y)\mathcal{H}(X,U)Z+\eta(Z)\mathcal{H}(X,Y)U\big{]}.
\eea
From (\ref{eq:4.4}) and (\ref{eq:4.7}) it follows that
\bea\nonumber
 (R(U,V)\cdot R)(X,Y)Z&=&\big{[}\{(\nabla_UR)(X,Y,Z,\phi V) - (\nabla_VR)(X,Y,Z,\phi U)\}\\
\label{eq:4.8} &+& \eta(X)H(Y,Z,U,V) -\eta(Y)H(X,Z,U,V)-2\eta(Z)H(X,Y,U,V)\\
\nonumber&+& \eta(U)H(X,Y,Z,V) - \eta(V)H(X,Y,Z,U)\big{]}\xi\\
\nonumber&+& \eta(U)\big{[}\eta(X)\mathcal{H}(V,Y)Z + \eta(Y)\mathcal{H}(X,V)Z +\eta(Z)\mathcal{H}(X,Y)V\big{]}\\
\nonumber&-& \eta(V)\big{[}\eta(X)\mathcal{H}(U,Y)Z + \eta(Y)\mathcal{H}(X,U)Z+\eta(Z)\mathcal{H}(X,Y)U\big{]}.
\eea
Thus in a locally $\phi$-semisymmetric Sasakian manifold, the relation (\ref{eq:4.8}) holds for any arbitrary vector fields $X, Y, Z, U$ and $V$ on $M$. Next, if the relation  (\ref{eq:4.8}) holds in a Sasakian manifold, then for any horizontal vector fields $X, Y, Z, U$ and $V$ on $M$, we get the relation (\ref{eq:3.18}) and hence the manifold is locally $\phi$-semisymmetric. Thus we can state the following:
\begin{theo}
A Sasakian manifold $M$ is locally $\phi$-semisymmetric if and only if the relation \textnormal{(\ref{eq:4.8})} holds for any arbitrary vector fields $X, Y, Z, U$ and $V$ on $M$.
\end{theo}
\begin{cor}\cite{tanno}
A semisymmetric Sasakian manifold is a manifold of constant curvature 1.
\end{cor}
\section{\bf Locally Ricci (resp., Projectively,  Conformally) $\phi$-Semisymmetric Sasakian Manifolds}
\begin{defi}
A Sasakian manifold $M$ is said to be a locally Ricci $\phi$-semisymmetric if the relation
\be\label{eq:5.1}
\phi^{2}[(R(U,V)\cdot Q)(X)]=0
\ee
holds for all horizontal vector fields $X, Y, Z, U$ and $V$ on $M$, $Q$ being the Ricci operator of the manifold.
\end{defi}
We know that
\be\label{eq:5.2}
(R(U,V)\cdot Q)(X)=R(U,V)QX-QR(U,V)X.
\ee
Applying $\phi^{2}$ on both sides of (\ref{eq:5.2}) we get
\be\label{eq:5.3}
\phi^{2}[(R(U,V)\cdot Q)(X)]=-(R(U,V)\cdot Q)(X)
\ee
for all horizontal vector fields $U, V$ and $X$ on $M$. This leads to the following:
\begin{theo}
A Sasakian manifold $M$ is locally Ricci $\phi$-semisymmetric if and only if  $(R(U,V)\cdot Q)(X)=0$ for all horizontal vector fields $U$, $V$ and $X$ on $M$.
\end{theo}
Now let $M$ be a locally  $\phi$-semisymmetric Sasakian manifold. Then the relation (\ref{eq:3.18}) holds on $M$. Taking inner product of  (\ref{eq:3.18}) with a horizontal vector field $W$ and then contracting over $X$ and $W$, we get $(R(U,V)\cdot S)(Y,Z)=0$ from which it follows that $(R(U,V)\cdot Q)(Y)=0$ for all horizontal vector fields $U$, $V$ and $Y$ on $M$. Thus in view of the Theorem 5.1, we can state the following:
\begin{theo}
A locally $\phi$-semisymmetric Sasakian manifold $M$ is locally Ricci $\phi$-semisymmetric.
\end{theo}
Now let  $U$, $V$ and $X$ are arbitrary vector fields on a Sasakian manifold $M$. Then in view of (\ref{eq:2.1}), (\ref{eq:2.4}), (\ref{eq:2.5}) and (\ref{eq:2.10}),  (\ref{eq:5.2}) yields
\bea\label{eq:5.4}
(R(\phi^2 U,\phi^2 V)\cdot Q)(\phi^2 X)&=&-(R(U,V)\cdot Q)(X)+\{E(X,V)\eta(U)-E(X,U)\eta(V)\}\xi\\
\nonumber&-&\eta(X)\{\eta(V)\mathcal{E}U-\eta(U)\mathcal{E}V\}
\eea
where $g(\mathcal{E}X,Y)=E(X,Y)$ and $E$ is given by 
\be\label{eq:5.5}
E(X,Y)=S(X,Y)-(n-1)g(X,Y).
\ee
Since $\phi^2 U$, $\phi^2 V$ and $\phi^2 X$ are orthogonal to $\xi$, in a locally Ricci $\phi$-semisymmetric Sasakian manifold $M$, from (\ref{eq:5.4}) we have 
\bea\label{eq:5.6}
(R(U,V)\cdot Q)(X)&=&\{E(X,V)\eta(U)-E(X,U)\eta(V)\}\xi\\
\nonumber&-&\eta(X)\{\eta(V)\mathcal{E}U-\eta(U)\mathcal{E}V\}.
\eea
Thus in a locally Ricci $\phi$-semisymmetric Sasakian manifold $M$ the relation (\ref{eq:5.6}) holds for any arbitrary vector fields $U$, $V$ and $X$ on $M$. Next, if the relation (\ref{eq:5.6}) holds in a Sasakian manifold $M$, then for all horizontal vector fields $U$, $V$ and $X$, we have $(R(U,V)\cdot Q)(X)=0$ and hence $M$ is locally Ricci $\phi$-semisymmetric. Thus we can state the following:
\begin{theo}
A Sasakian manifold $M$ is locally Ricci $\phi$-semisymmetric if and only if  the relation \textnormal{(\ref{eq:5.6})} holds for any arbitrary vector fields $U$, $V$ and $X$ on $M$.
\end{theo}
\begin{cor}\cite{tanno}
A  Ricci semisymmetric Sasakian manifold is an Einstein manifold.
\end{cor}
\begin{defi}
A Sasakian manifold $M$ is said to be a locally projectively (resp. conformally) $\phi$-semisymmetric if the relation
\be\label{eq:5.7}
\phi^{2}[(R(U,V)\cdot P)(X,Y)Z] \ \big{(}\textnormal{resp.} \  \phi^{2}[(R(U,V)\cdot C)(X,Y)Z]\big{)}=0
\ee
holds for all horizontal vector fields $X, Y, Z, U$ and $V$ on $M$, $P$ (resp. $C$) being the projective (resp. conformal) curvature tensor of the manifold.
\end{defi}
The projective transformation is such that geodesics transformed into geodesics \cite{weyl} and as the invariant of such transformation the Weyl projective curvature tensor $P$ of type (1,3) is given by \cite{weyl}
\be
\label{eq:5.8}
P(X,Y)Z = R(X,Y)Z - \frac{1}{n-1} \big{[}S(Y,Z)X-S(X,Z)Y\big{]}.
\ee
The conformal transformation is an angle preserving mapping and as the invariant of such transformation the Weyl conformal curvature tensor $C$ of type (1,3) on a Riemannian manifold $M$, $n>3$, is given by \cite{weyl}
\bea
\label{eq:5.9}
C(X,Y)Z &=& R(X,Y)Z - \frac{1}{n-2} \{S(Y,Z)X-S(X,Z)Y+g(Y,Z)QX\\
\nonumber&-&g(X,Z)QY\}+\frac{r}{(n-1)(n-2)}\{g(Y,Z)X-g(X,Z)Y\}.
\eea
From (\ref{eq:5.8}) we get 
\bea\label{eq:5.10} 
\ \ \ \ \ \  (R(U,V)\cdot P)(X,Y)Z&&=(R(U,V)\cdot R)(X,Y)Z-\frac{1}{n-1}\left[(R(U,V)\cdot S)(Y,Z)X\right.\\
\nonumber&&\left.-(R(U,V)\cdot S)(X,Z)Y\right].
\eea
Applying $\phi^{2}$ on both sides of (\ref{eq:5.10}) and using (\ref{eq:3.9}) we obtain
\bea\label{eq:5.11}\ \ \ \ \ \ \ 
\phi^{2}[(R(U,V)\cdot P)(X,Y)Z]&&=-(R(U,V)\cdot R)(X,Y)Z\\
\nonumber&&+\left[(\nabla_UR)(X,Y,Z,\phi V) - (\nabla_VR)(X,Y,Z,\phi U)\right]\xi \\
\nonumber&&+\frac{1}{n-1}\left[(R(U,V)\cdot S)(Y,Z)X-(R(U,V)\cdot S)(X,Z)Y\right]
\eea
for all horizontal vector fields $X, Y, Z, U$ and $V$ on $M$.\\
\indent Now we suppose that $M$ is a locally projectively $\phi$-semisymmetric Sasakian manifold. Then from (\ref{eq:5.11}) we obtain
\bea\label{eq:5.12}
(R(U,V)\cdot R)(X,Y)Z&&=\left[(\nabla_UR)(X,Y,Z,\phi V) - (\nabla_VR)(X,Y,Z,\phi U)\right]\xi \\
\nonumber&&+\frac{1}{n-1}\left[(R(U,V)\cdot S)(Y,Z)X-(R(U,V)\cdot S)(X,Z)Y\right]
\eea
for all horizontal vector fields $X, Y, Z, U$ and $V$ on $M$. Taking inner product of  (\ref{eq:5.12})  with a horizontal vector field $W$ and then contracting over $X$ and $Z$, we get 
\be\label{eq:5.13}
(R(U,V)\cdot S)(Y,W)=0
\ee
 for all horizontal vetor fields $U, V, Y$ and $W$ on $M$ and hence by Theorem 5.1,  it follows that the manifold $M$ is locally Ricci $\phi$-semisymmetric. Using (\ref{eq:5.13}) in (\ref{eq:5.12}), it follows that the manifold $M$ is locally $\phi$-semisymmetric.\\
\indent Next, we suppose that $M$ is a locally $\phi$-semisymmetric Sasakian manifold. Then the relation (\ref{eq:3.18}) holds on $M$. Taking inner product of  (\ref{eq:3.18}) with a horizontal vector field $W$ and then contracting over $X$ and $W$, we get $(R(U,V)\cdot S)(Y,Z)=0$ for all horizontal vetor fields $U, V, Y$ and $Z$ on $M$ and hence from (\ref{eq:5.11}) it follows that the manifold $M$ is locally projectively $\phi$-semisymmetric. This leads to the following:
\begin{theo}
A locally projectively $\phi$-semisymmetric Sasakian manifold $M$ is locally $\phi$-semi-symmetric and vice versa.
\end{theo}
Now from (\ref{eq:5.9}) we get 
\bea\label{eq:5.14} 
(R(U,V)\cdot C)(X,Y)Z&=&(R(U,V)\cdot R)(X,Y)Z\\
\nonumber&-&\frac{1}{n-2}[(R(U,V)\cdot S)(Y,Z)X-(R(U,V)\cdot S)(X,Z)Y\\
\nonumber&+&g(Y,Z)(R(U,V)\cdot Q)(X)-g(X,Z)(R(U,V)\cdot Q)(Y)].
\eea
Applying $\phi^{2}$ on both sides of (\ref{eq:5.14}) and using (\ref{eq:3.9}) and (\ref{eq:5.3}) we obtain
\bea\label{eq:5.15}
\ \ \ \ \phi^{2}[(R(U,V)\cdot C)(X,Y)Z]&=&-(R(U,V)\cdot R)(X,Y)Z\\
\nonumber&+&[(\nabla_UR)(X,Y,Z,\phi V) - (\nabla_VR)(X,Y,Z,\phi U)]\xi \\
\nonumber&+&\frac{1}{n-2}\left[(R(U,V)\cdot S)(Y,Z)X-(R(U,V)\cdot S)(X,Z)Y\right.\\
\nonumber&+&\left.g(Y,Z)(R(U,V)\cdot Q)(X)-g(X,Z)(R(U,V)\cdot Q)(Y)\right]
\eea
for all horizontal vector fields $X, Y, Z, U$ and $V$ on $M$. This leads to the following:
\begin{theo}
A Sasakian manifold $M$ is locally conformally $\phi$-semisymmetric if and only if the relation 
\bea\label{eq:5.16}
(R(U,V)\cdot R)(X,Y)Z&=&\left[(\nabla_UR)(X,Y,Z,\phi V) - (\nabla_VR)(X,Y,Z,\phi U)\right]\xi \\
\nonumber&+&\frac{1}{n-2}\left[(R(U,V)\cdot S)(Y,Z)X-(R(U,V)\cdot S)(X,Z)Y\right.\\
\nonumber&+&\left.g(Y,Z)(R(U,V)\cdot Q)(X)-g(X,Z)(R(U,V)\cdot Q)(Y)\right]
\eea
holds for all horizontal vector fields $X, Y, Z, U$ and $V$ on $M$.
\end{theo}
Let $M$ be a locally  $\phi$-semisymmetric Sasakian manifold. Then $M$ is locally Ricci $\phi$-semi-symmetric and thus in view of  (\ref{eq:3.18}), it follows from (\ref{eq:5.15}) that $\phi^{2}[(R(U,V)\cdot C)(X,Y)Z]=0$ for all horizontal vector fields $X, Y, Z, U$ and $V$ on $M$. Hence the manifold $M$ is locally conformally $\phi$-semisymmetric.\\
\indent Again, we consider  $M$ as the locally conformally  $\phi$-semisymmetric Sasakian manifold. If $M$ is locally  $\phi$-semisymmetric Sasakian manifold, then from (\ref{eq:5.16}) it follows that $(R(U,V)\cdot S)(Y,Z)=0$ which implies that $(R(U,V)\cdot Q)(Y)=0$ for all horizontal vector fields $U$, $V$ and $Y$ on $M$ and hence by Theorem 5.1, the manifold $M$ is locally Ricci $\phi$-semisymmetric. Again, if $M$ is locally Ricci $\phi$-semisymmetric, then $(R(U,V)\cdot Q)(Y)=0$ for all horizontal vector fields $U$, $V$ and $Y$ on $M$ and hence by Theorem 3.3, it follows from (\ref{eq:5.16}) that the manifold $M$ is locally  $\phi$-semisymmetric. This leads to the following:
\begin{theo}
A locally $\phi$-semisymmetric Sasakian manifold $M$ is locally conformally $\phi$-semisymmetric. The converse is true if and only if the manifold $M$ is locally Ricci $\phi$-semi-symmetric.
\end{theo}
Now let $X, Y, Z, U$ and $V$ be any arbitrary vector fields on a Sasakian manifold $M$. Then using (\ref{eq:2.1}), (\ref{eq:2.10}),  (\ref{eq:4.2}) and (\ref{eq:5.16}) we obtain 
\bea\nonumber
(R(\phi^{2}U,\phi^{2}V)\cdot R)(\phi^{2}X,\phi^{2}Y)\phi^{2}Z&=&[\{(\nabla_VR)(X,Y,Z,\phi U)-(\nabla_UR)(X,Y,Z,\phi V)\}\\
\nonumber&-&\eta(X)H(Y,Z,U,V)+\eta(Y)H(X,Z,U,V)\\
\label{eq:5.17}&+& 2\eta(Z)H(X,Y,U,V)]\xi\\
\nonumber&-&\frac{1}{n-2}\big{[}\{(R(U,V)\cdot S)(Y,Z)X-(R(U,V)\cdot S)(X,Z)Y\}\\
\nonumber&-&\big{\{}(R(U,V)\cdot S)(Y,Z)\eta(X)-(R(U,V)\cdot S)(X,Z)\eta(Y)\big{\}}\xi\\
\nonumber&-&\big{\{}E(V,Z)\eta(U)-E(U,Z)\eta(V)\big{\}}\{\eta(Y)X-\eta(X)Y\}\\
\nonumber&+&\big{\{}E(Y,U)X-E(X,U)Y\big{\}}\eta(Z)\eta(V)\\
\nonumber&-&\big{\{}E(Y,V)X-E(X,V)Y\big{\}}\eta(Z)\eta(U)\big{]}-\frac{1}{n-2}\\
\nonumber&&\big{[}\{g(Y,Z)(R(U,V).Q)(X)-g(X,Z)(R(U,V).Q)(Y)\}\\
\nonumber&-&\{\eta(Y)(R(U,V).Q)(X)-\eta(X)(R(U,V).Q)(Y)\}\eta(Z)\\
\nonumber&+&\{g(Y,Z)\eta(X)-g(X,Z)\eta(Y)\}\{\eta(V)\mathcal{E}U-\eta(U)\mathcal{E}V\}\\
\nonumber&-&\big{\{}E(X,V)\eta(U)-E(X,U)\eta(V)\big{\}}g(Y,Z)\xi\\
\nonumber&+&\big{\{}E(Y,V)\eta(U)-E(Y,U)\eta(V)\big{\}}g(X,Z)\xi\big{]}
\eea
where $g(\mathcal{E}U,V)=E(U,V)$ and $E$ is given by (\ref{eq:5.5}).\\
From (\ref{eq:5.17}) and (\ref{eq:4.7}) it follows that 
\bea\label{eq:5.18} 
\ \ \ \ \ \ \ (R(U,V)\cdot R)(X,Y)Z&=&\big{[}\{(\nabla_UR)(X,Y,Z,\phi V) - (\nabla_VR)(X,Y,Z,\phi U) \}\\
\nonumber&+& \eta(X)H(Y,Z,U,V)- \eta(Y)H(X,Z,U,V)\\
\nonumber&-& 2\eta(Z)H(X,Y,U,V)+ \eta(U)H(X,Y,Z,V)\\
\nonumber&-& \eta(V)H(X,Y,Z,U)\big{]}\xi+ \eta(U)\big{[}\eta(X)\mathcal{H}(V,Y)Z\\
\nonumber&+& \eta(Y)\mathcal{H}(X,V)Z+\eta(Z)\mathcal{H}(X,Y)V\big{]}\\ 
\nonumber&-& \eta(V)\big{[}\eta(X)\mathcal{H}(U,Y)Z+ \eta(Y)\mathcal{H}(X,U)Z+\eta(Z)\mathcal{H}(X,Y)U\big{]}\\
\nonumber&+& \frac{1}{n-2}\big{[}\{(R(U,V)\cdot S)(Y,Z)X-(R(U,V)\cdot S)(X,Z)Y\}\\
\nonumber&-&\big{\{}(R(U,V)\cdot S)(Y,Z)\eta(X)-(R(U,V)\cdot S)(X,Z)\eta(Y)\big{\}}\xi\\
\nonumber&-& \big{\{}E(V,Z)\eta(U)-E(U,Z)\eta(V)\big{\}}\{\eta(Y)X-\eta(X)Y\}\\
\nonumber&+& \big{\{}E(Y,U)X-E(X,U)Y\big{\}}\eta(Z)\eta(V)
\eea
\bea
\nonumber&-& \big{\{}E(Y,V)X-E(X,V)Y\big{\}}\eta(Z)\eta(U)\big{]}\\
\nonumber&+& \frac{1}{n-2}\big{[}\{g(Y,Z)(R(U,V).Q)(X)-g(X,Z)(R(U,V).Q)(Y)\}\\
\nonumber&-& \{\eta(Y)(R(U,V).Q)(X)-\eta(X)(R(U,V).Q)(Y)\}\eta(Z)\\
\nonumber&+& \{g(Y,Z)\eta(X)-g(X,Z)\eta(Y)\}\{\eta(V)\mathcal{E}U-\eta(U)\mathcal{E}V\}\\
\nonumber&-& \big{\{}E(X,V)\eta(U)-E(X,U)\eta(V)\big{\}}g(Y,Z)\xi\\
\nonumber&+& \big{\{}E(Y,V)\eta(U)-E(Y,U)\eta(V)\big{\}}g(X,Z)\xi\big{]}
\eea
where $H(X,Y,Z,U)=g(\mathcal{H}(X,Y)Z,U)$ and $g(\mathcal{E}U,V)=E(U,V)$, $\mathcal{H}$ and $E$ are given by (\ref{eq:4.5}) and (\ref{eq:5.5}) respectively.
Thus in a locally conformally $\phi$-semisymmetric Sasakian manifold $M$ the relation (\ref{eq:5.18}) holds for any arbitrary vector fields $X, Y, Z, U$ and $V$ on $M$. Next, if the relation (\ref{eq:5.18}) holds in a Sasakian manifold $M$, then for all horizontal vector fields $X, Y, Z, U$ and $V$ on $M$, we have (\ref{eq:5.16}), that is, the manifold is locally conformally $\phi$-semisymmetric. This leads to the following:
\begin{theo}
A Sasakian manifold $M$ is locally conformally $\phi$-semisymmetric if and only if the relation \textnormal{(\ref{eq:5.18})} holds for any arbitrary vector fields $X, Y, Z, U$ and $V$ on $M$.   
\end{theo}
\begin{cor}\cite{chaki}
A conformally semisymmetric Sasakian manifold is a manifold of constant curvature 1.
\end{cor}
\begin{rem}
Since the skew-symmetric operator $R(X,Y)$ and the structure tensor $\phi$ of the Sasakian manifold both are commutes with the contraction, it follows from Theorem 6.6(ii) of Shaikh and Kundu \cite{kundu} that the same conclusion of the Theorem 5.5, Theorem 5.6 and Theorem 5.7 holds for locally conharmonically $\phi$-semisymmetric Sasakian manifold.
\end{rem}
Again, by linear combination of $R$, $S$ and $g$, Shaikh and Kundu \cite{kundu} defined a generalized curvature tensor $B$ (see, equation (2.1) of \cite{kundu}) of type (1,3), called $B$-tensor which includes various curvature tensors as particular cases. Then Shaikh and Kundu (see, equation (5.5) of \cite{kundu}) shows that this $B$-tensor turns into the following form:
\bea\label{eq:5.19}
B(X,Y)Z&=&b_{0}R(X,Y)Z+b_{1}\{S(Y,Z)X-S(X,Z)Y\\
\nonumber&+&g(Y,Z)QX-g(X,Z)QY\}+b_{2}r\{g(Y,Z)X-g(X,Z)Y\}
\eea
where $b_{0}$, $b_{1}$ and $b_{2}$ are scalars.\\
We note that if\\
(a)  $b_{0}=1$, $b_{1}=0$ and $b_{2}=-\frac{1}{n(n-1)}$;\\
(b) $b_{0}=1$, $b_{1}=-\frac{1}{(n-2)}$ and $b_{2}=\frac{1}{(n-1)(n-2)}$;\\
(c)  $b_{0}=1$, $b_{1}= -\frac{1}{(n-2)}$ and $b_{2}=0$;\\ and
(d)  $b_{2}=-\frac{1}{n}\big{(}\frac{b_{0}}{n-1}+2b_{1}\big{)}$,\\
then from (\ref{eq:5.19}) it follows that  the  $B$-tensor turns into the (a) concircular, (b) conformal, (c) conharmonic and (d) quasi-conformal curvature tensor respectively. For details about the $B$-tensor we refer the reader to see Shaikh and Kundu \cite{kundu} and also references therein.
\begin{defi}
A Sasakian manifold $M$ is said to be a locally $B$-$\phi$-semisymmetric if the relation
\be\label{eq:5.20}
\phi^{2}[(R(U,V)\cdot B)(X,Y)Z]=0
\ee
holds for all horizontal vector fields $X, Y, Z, U$ and $V$ on $M$, $B$  being the generalized curvature tensor of the manifold.
\end{defi}
From (\ref{eq:5.19}) we get 
\bea\label{eq:5.21} 
\ \ \ \ \ \  (R(U,V)\cdot B)(X,Y)Z&&=b_{0}(R(U,V)\cdot R)(X,Y)Z+b_{1}\left[(R(U,V)\cdot S)(Y,Z)X\right.\\
\nonumber&&\left.-(R(U,V)\cdot S)(X,Z)Y+g(Y,Z)(R(U,V)\cdot Q)(X)\right.\\
\nonumber&&\left.-g(X,Z)(R(U,V)\cdot Q)(Y)\right].
\eea
Applying $\phi^{2}$ on both sides of (\ref{eq:5.21}) and using (\ref{eq:3.9}) and (\ref{eq:5.3}) we obtain
\bea\label{eq:5.22}\ \ \ \ \ \ 
\phi^{2}[(R(U,V)\cdot B)(X,Y)Z]&=&-b_{0}[(R(U,V)\cdot R)(X,Y)Z\\
\nonumber&-&\left\{(\nabla_UR)(X,Y,Z,\phi V) - (\nabla_VR)(X,Y,Z,\phi U)\right\}\xi] \\
\nonumber&-&b_{1}\left[(R(U,V)\cdot S)(Y,Z)X-(R(U,V)\cdot S)(X,Z)Y\right.\\
\nonumber&+&\left.g(Y,Z)(R(U,V)\cdot Q)(X)-g(X,Z)(R(U,V)\cdot Q)(Y)\right]
\eea
for all horizontal vector fields $X, Y, Z, U$ and $V$ on $M$. This leads to the following:
\begin{theo}
A Sasakian manifold $M$ is locally $B$-$\phi$-semisymmetric if and only if  
\bea\label{eq:5.23}
(R(U,V)\cdot R)(X,Y)Z&=&\left[(\nabla_UR)(X,Y,Z,\phi V) - (\nabla_VR)(X,Y,Z,\phi U)\right]\xi \\
\nonumber&-&\frac{b_{1}}{b_{0}}\left[(R(U,V)\cdot S)(Y,Z)X-(R(U,V)\cdot S)(X,Z)Y\right.\\
\nonumber&&\left.+g(Y,Z)(R(U,V)\cdot Q)(X)-g(X,Z)(R(U,V)\cdot Q)(Y)\right]
\eea
for all horizontal vector fields $X, Y, Z, U$ and $V$ on $M$, provided $b_{0}\neq 0$.
\end{theo}
Now taking inner product of  (\ref{eq:5.23}) with a horizontal vector field $W$ and then contracting over $X$ and $W$, we get 
\be\label{eq:5.24} 
\{b_{0}+(n-2)b_{1}\}(R(U,V)\cdot S)(Y,Z)=0
\ee
for all horizontal vector fields $X, Y, Z, U$ and $V$ on $M$.\\
From (\ref{eq:5.24}) following two cases arise:\\
\textbf{Case-I:} If $b_{0}+(n-2)b_{1}\neq 0$, then from (\ref{eq:5.24}) we have 
\be\label{eq:5.25} 
(R(U,V)\cdot S)(Y,Z)=0
\ee
for all horizontal vector fields $X, Y, Z, U$ and $V$ on $M$, from which it follows that $(R(U,V)\cdot Q)(Y)=0$ for all horizontal vector fields $U$, $V$ and $Y$ on $M$. This leads to the following:
\begin{theo}
A locally  $B$-$\phi$-semisymmetric Sasakian manifold $M$ is locally Ricci $\phi$-semi-symmetric provided that $b_{0}+(n-2)b_{1}\neq 0$.
\end{theo}
\begin{cor}
A locally concircularly  $\phi$-semisymmetric Sasakian manifold $M$ is locally Ricci $\phi$-semisymmetric.
\end{cor}
\begin{cor}
A locally quasi-conformally $\phi$-semisymmetric Sasakian manifold $M$ is locally Ricci $\phi$-semisymmetric  provided that $b_{0}+(n-2)b_{1}\neq 0$.
\end{cor}
Now if $b_{0}+(n-2)b_{1}\neq 0$, then in view of  (\ref{eq:5.25}), (\ref{eq:5.23}) takes the form (\ref{eq:3.18})
for all horizontal vector fields $X, Y, Z, U$ and $V$ on $M$ and hence the manifold $M$ is locally  $\phi$-semisymmetric. Again, if we consider the manifold $M$ as locally $\phi$-semisymmetric, then the relation (\ref{eq:3.18}) holds on $M$. Taking inner product of  (\ref{eq:3.18}) with a horizontal vector field $W$ and then contracting over $X$ and $W$, we get $(R(U,V)\cdot S)(Y,Z)=0$ for all horizontal vetor fields $U, V, Y$ and $Z$ on $M$ and hence from (\ref{eq:5.22}) it follows that the manifold $M$ is locally $B$-$\phi$-semisymmetric. Thus we can state the following:
\begin{theo}
In a Sasakian manifold $M$, local  $B$-$\phi$-semisymmetry  and local  $\phi$-semisymmetry are equivalent provided that $b_{0}+(n-2)b_{1}\neq 0$.
\end{theo}
\begin{cor}
In a Sasakian manifold $M$, local concircular  $\phi$-semisymmetry  and local  $\phi$-semisymmetry are equivalent.
\end{cor}
\begin{cor}
In a Sasakian manifold $M$, local quasi-conformal $\phi$-semisymmetry  and local  $\phi$-semisymmetry are equivalent  provided that $b_{0}+(n-2)b_{1}\neq 0$.
\end{cor}
\begin{rem}
Since the skew-symmetric operator $R(X,Y)$ and the structure tensor $\phi$ of the Sasakian manifold both are commutes with the contraction, it follows from Theorem 6.6(i) of Shaikh and Kundu \cite{kundu} that the same conclusion of the Corollary 5.3 and Corollary 5.5 holds for locally projectively $\phi$-semisymmetric Sasakian manifold as the contraction on projective curvature tensor gives rise the Ricci operator although projective curvature tensor is not a generalized curvature tensor.
\end{rem}
\noindent \textbf{Case-II:} If $b_{0}+(n-2)b_{1}= 0$, then from (\ref{eq:5.22}) we have
\bea\label{eq:5.26}\ \ \ \ \ \ 
\phi^{2}[(R(U,V)\cdot B)(X,Y)Z]&=&-b_{0}[(R(U,V)\cdot R)(X,Y)Z\\
\nonumber&-&\left\{(\nabla_UR)(X,Y,Z,\phi V) - (\nabla_VR)(X,Y,Z,\phi U)\right\}\xi] \\
\nonumber&+&\frac{b_{0}}{n-2}\left[(R(U,V)\cdot S)(Y,Z)X-(R(U,V)\cdot S)(X,Z)Y\right.\\
\nonumber&+&\left.g(Y,Z)(R(U,V)\cdot Q)(X)-g(X,Z)(R(U,V)\cdot Q)(Y)\right]
\eea
for all horizontal vector fields $X, Y, Z, U$ and $V$ on $M$. This leads to the following:
\begin{theo}
A Sasakian manifold $M$ is locally $B$-$\phi$-semisymmetric if and only if  
\bea\label{eq:5.27}
(R(U,V)\cdot R)(X,Y)Z&=&\left[(\nabla_UR)(X,Y,Z,\phi V) - (\nabla_VR)(X,Y,Z,\phi U)\right]\xi \\
\nonumber&+&\frac{1}{n-2}\left[(R(U,V)\cdot S)(Y,Z)X-(R(U,V)\cdot S)(X,Z)Y\right.\\
\nonumber&&\left.+g(Y,Z)(R(U,V)\cdot Q)(X)-g(X,Z)(R(U,V)\cdot Q)(Y)\right]
\eea
for all horizontal vector fields $X, Y, Z, U$ and $V$ on $M$ provided that $b_{0}+(n-2)b_{1}= 0$.
\end{theo}
\begin{cor}
A Sasakian manifold $M$ is locally conformally (resp. conharmonically) $\phi$-semisymmetric if and only if the relation (\ref{eq:5.27}) holds.
\end{cor}
\begin{cor}
A Sasakian manifold $M$ is locally quasi-conformally  $\phi$-semisymmetric if and only if the relation (\ref{eq:5.27}) holds provided that $b_{0}+(n-2)b_{1}= 0$.
\end{cor}
Let $M$ be a locally  $\phi$-semisymmetric Sasakian manifold. Then $M$ is locally Ricci $\phi$-semi-symmetric and thus in view of  (\ref{eq:3.18}), it follows from (\ref{eq:5.22}) that $\phi^{2}[(R(U,V)\cdot B)(X,Y)Z]=0$ for all horizontal vector fields $X, Y, Z, U$ and $V$ on $M$. Hence the manifold $M$ is locally $B$-$\phi$-semisymmetric.\\
\indent Again, we consider $M$ as the locally $B$-$\phi$-semisymmetric Sasakian manifold. If $b_{0}+(n-2)b_{1}\neq 0$, then $M$ is locally $\phi$-semisymmetric. So we suppose that $b_{0}+(n-2)b_{1}= 0$. If  $M$ is locally  $\phi$-semisymmetric, then from (\ref{eq:5.27}) it follows that $(R(U,V)\cdot S)(Y,Z)=0$, which implies that $(R(U,V)\cdot Q)(Y)=0$ for all horizontal vector fields $U$, $V$ and $Y$ on $M$. Thus in view of  Theorem 5.1,  the manifold $M$ is locally Ricci $\phi$-semisymmetric. Again, if $M$ is locally Ricci $\phi$-semisymmetric, then $(R(U,V)\cdot Q)(Y)=0$ for all horizontal vector fields $U$, $V$ and $Y$ on $M$. Thus in view of Theorem 3.3, it follows from (\ref{eq:5.27}) that the manifold $M$ is locally  $\phi$-semisymmetric. This leads to the following:
\begin{theo}
A locally $\phi$-semisymmetric Sasakian manifold $M$ is locally $B$-$\phi$-semisymmetric. The converse is true for $b_{0}+(n-2)b_{1}= 0$ if and only if the manifold $M$ is locally Ricci $\phi$-semisymmetric.
\end{theo}
If $X, Y, Z, U$ and $V$ are  arbitrary vector fields on  $M$, then proceeding similarly as in the case of conformal curvature tensor, it is easy to check that (\ref{eq:5.18}) holds for $b_{0}+(n-2)b_{1}= 0$. Hence we can state the following:
\begin{theo}
A Sasakian manifold $M$ is locally $B$-$\phi$-semisymmetric if and only if the relation \textnormal{(\ref{eq:5.18})} holds for any arbitrary vector fields $X, Y, Z, U$ and $V$ on $M$ provided that $b_{0}+(n-2)b_{1}= 0$.
\end{theo}
\begin{cor}
A Sasakian manifold $M$ is locally conformally (resp. conharmonically) $\phi$-semisymmetric if and only if the relation \textnormal{(\ref{eq:5.18})} holds for any arbitrary vector fields $X, Y, Z, U$ and $V$ on $M$.
\end{cor}
\begin{cor}
A Sasakian manifold $M$ is locally quasi-conformally $\phi$-semisymmetric if and only if the relation \textnormal{(\ref{eq:5.18})} holds for any arbitrary vector fields $X, Y, Z, U$ and $V$ on $M$ provided that $b_{0}+(n-2)b_{1}= 0$.
\end{cor}
\textbf{Conclusion.}  From the above discussion and results, we conclude that the study of local $\phi$-semisymmetry is meaningful as a generalized notion of local $\phi$-symmetry and semisymmetry. From Theorem 6.6 and Corollary 6.2 of Shaikh and Kundu \cite{kundu} we also conclude that the same characterization of local $\phi$-semisymmetry of a Sasakian manifold holds for the locally projectively $\phi$-semisymmetric and locally concircularly $\phi$-semisymmetric Sasakian manifolds as the contraction on projective or concircular curvature tensor gives rise the Ricci operator. And also from Theorem 6.6 and Corollary 6.2 of Shaikh and Kundu \cite{kundu} we again conclude that the local conformal $\phi$-semisymmetry and local conharmonical $\phi$-semisymmetry on a  Sasakian manifold are equivalent. However, the study of local $\phi$-semisymmetry and local conformal $\phi$-semisymmetry  are meaningful as they are not equivalent. Finally, we conclude that the study of local $\phi$-semisymmetry on a Sasakian manifold by considering any other generalized curvature tensor of type (1,3)(which are the linear combination of $R$, $S$ and $g$ ) is either meaningless or redundant due to their equivalency.

 \end{document}